\documentclass{article}
\usepackage{amsmath}
\usepackage{harvard}
\usepackage{amsthm}
\usepackage{eucal}

\newcounter{ctheorem}[section]                   
\newtheorem{lemma}[ctheorem]{Lemma}              
\newtheorem{thm}[ctheorem]{Theorem}              
\newtheorem{defin}[ctheorem]{Definition}         
\theoremstyle{definition}

\makeatletter
\@addtoreset{equation}{section}
\makeatother

\def\roundkern#1#2#3#4{\copy255 \kern-#1\wd255 \vrule width #4\wd255
    height #2\ht255 depth #3\ht255 \kern#1\wd255}
\def\straightletter#1{{\mathord{\rm I\mkern-3.6mu #1}}}
\def\N{\straightletter{N}}
\def\R{\straightletter{R}}
\def\Z{\mathord{\rm Z\mkern-6.1mu Z}}
\newcommand{\order}{\CMcal{O}}                   
\newcommand{\nodes}{\CMcal{A}}                   
\newcommand{\To}{\rightarrow}                    
\newcommand{\dotz}[2]{#1,\dotsc,#2}              
\newcommand{\eps}{\varepsilon}                   
\newcommand{\h}[1]{\widehat{#1}}                 
\newcommand{\ehat}{\,\widehat{\;}\,}             
\DeclareMathOperator{\support}{supp}             
\newcommand{\supp}[1]{\support{#1}}              
\DeclareMathOperator{\volume}{vol}               
\newcommand{\vol}[1]{\volume(#1)}                
\newcommand{\set}[1]{\{#1\}}                     
\newcommand{\abs}[1]{\lvert#1\rvert}             
\newcommand{\norm}[1]{\lVert#1\rVert}            
\newcommand{\biggabs}[1]{\biggl\lvert#1\biggr\rvert}

\begin{document}


\global\baselineskip22pt plus2pt minus2pt

\begin{center}
    {\Large \bf Approximation Orders for Interpolation by Surface Splines to Rough Functions}
\end{center}
\vspace*{5mm}
\begin{center}
    Rob Brownlee\footnote{This author was supported by a studentship
    from the Engineering and Physical Sciences Research Council.} and Will Light\footnote{
    The first author would like to dedicate this paper to the memory of Will Light.} \\
    Department of Mathematics and Computer Science, University of
    Leicester, University Road, Leicester LE1 7RH, England.
\end{center}
\vspace*{5mm}
%
%

\begin{flushleft}
{\large \bf Abstract}
\end{flushleft}
In this paper we consider the approximation of functions by radial
basic function interpolants. There is a plethora of results about
the asymptotic behaviour of the error between appropriately smooth
functions and their interpolants, as the interpolation points fill
out a bounded domain in $\R^d$. In all of these cases, the analysis
takes place in a natural function space dictated by the choice of
radial basic function -- the native space. In many cases, the native
space contains functions possessing a certain amount of smoothness.
We address the question of what can be said about these error
estimates when the function being interpolated fails to have the
required smoothness. These are the rough functions of the title. We
limit our discussion to surface splines, as an exemplar of a wider
class of radial basic functions, because we feel our techniques are
most easily seen and understood in this setting.
\newpage
%
%

\section{INTRODUCTION}

\noindent The process of interpolation by translates of a basic
function is a popular tool for the reconstruction of a multivariate
function from a scattered data set. The setup of the problem is as
follows. We are supplied with a finite set of interpolation points
$\nodes \subset \R^d$ and a function $f:\nodes \To \R$. We wish to
construct an interpolant to $f$ of the form
\begin{equation}\label{interpolant form}
(Sf)(x) = \sum_{a \in \nodes} \mu_a \psi(x-a) + p(x),\qquad
\mbox{for $x \in \R^d$.}
\end{equation}
Here, $\psi$ is a real-valued function defined on $\R^d$, and the
principle ingredient of our interpolant is the use of the translates
of $\psi$ by the points in $\nodes$. The function $\psi$ is referred
to as the {\it basic function}. The function $p$ in
Equation~\eqref{interpolant form} is a polynomial on $\R^d$ of total
degree at most $k-1$. The linear space of all such polynomials will
be denoted by $\Pi_{k-1}$. Of course, for $Sf$ to interpolate $f$
the real numbers $\mu_a$ and the polynomial $p$ must be chosen to
satisfy the system
\begin{equation*}
(Sf)(a) = f(a),\qquad \mbox{for $a \in \nodes$}.
\end{equation*}
It is natural to desire a unique solution to the above system.
However, with the present setup, there are less conditions available
to determine $Sf$ than there are free parameters in $Sf$. There is a
standard way of determining the remaining conditions, which are
often called the {\it natural boundary conditions}:
$$
\sum_{a\in\nodes} \mu_a q(a) = 0, \qquad \mbox{for all $q\in
\Pi_{k-1}$}.
$$
It is now essential that $\nodes$ is $\Pi_{k-1}$--unisolvent. This
means that if $q \in \Pi_{k-1}$ vanishes on $\nodes$ then $q$ must
be zero. Otherwise the polynomial term can be adjusted by any
polynomial which is zero on $\nodes$. However, more conditions are
needed to ensure uniqueness of the interpolant. The requirement that
$\psi$ should be strictly conditionally positive definite of order
$k$ is one possible assumption. To see explanations of why these
conditions arise, the reader is directed
to~\citeasnoun{cheneylight}. In most of the common applications the
function $\psi$ is a radial function. That is, there is a function
$\phi:\R_+ \rightarrow \R$ such that $\psi = \phi \circ
\abs{\,\cdot\,}$, where $\abs{\,\cdot\,}$ is the Euclidean norm. In
these cases we refer to $\psi$ as a {\it radial basic function}.

\citename{Duchon} \citeyear{Duchon,Duchon2} was amongst the first to
study interpolation problems of this flavour. His approach was to
formulate the interpolation problem as a variational one. To do this
we assume we have a space of continuous functions $X$ which carries
a seminorm $\abs{\,\cdot\,}$. The so-called minimal norm interpolant
to $f \in X$ on $\nodes$ from $X$ is the function $Sf \in X$
satisfying
\begin{enumerate}
\item $(Sf)(a)=f(a)$, for all $a \in \nodes$;

\item $\abs{Sf} \leq \abs{g}$, for all $g \in X$ such that
$g(a)=f(a)$ for all $a \in \nodes$.
\end{enumerate}
The spaces that Duchon considers are in fact spaces of tempered
distributions which he is able to embed in $C(\R^d)$. Let ${\cal
S}'$ be the space of all tempered distributions on $\R^d$. The
particular spaces of distributions that we will be concerned with
are called Beppo-Levi spaces. The $k$\textsuperscript{th} order
Beppo-Levi space is denoted by $BL^k(\Omega)$ and defined as
\begin{equation*}
BL^k(\Omega) = \left\{ f\in{\cal S}' : \mbox{$D^\alpha f \in
L_2(\Omega)$, $\alpha \in \Z^d_+$, $\abs{\alpha}=k$} \right\},
\end{equation*}
with seminorm
\begin{equation*}
    \abs{f}_{k,\Omega}=\Biggl( \sum_{\abs{\alpha}=k} c_\alpha
    \int_{\Omega}\abs{(D^\alpha f)(x)}^2\, dx\Biggr)^{1/2},\qquad f\in
BL^k(\Omega).
\end{equation*}
The constants $c_\alpha$ are chosen so that the seminorm is
rotationally invariant:
\begin{equation*}
    \sum_{\abs{\alpha}=k} c_\alpha x^{2\alpha} = \abs{x}^{2k}, \qquad \mbox{for
all $x\in\R^d$}.
\end{equation*}
We assume throughout the paper that $2k>d$, because this has the
affect that $BL^k(\Omega)$ is embedded in the continuous
functions~\cite{Duchon}. The spaces $BL^k(\R^d)$ give rise to
minimal norm interpolants which are exactly of the form given in
Equation~\eqref{interpolant form}, where the radial basic function
is $x\mapsto\abs{x}^{2k-d}$ or $x\mapsto \abs{x}^{2k-d}\
\log{\abs{x}}$, depending on the parity of $d$.

It is perhaps no surprise to learn that the related functions $\psi$
are strictly conditionally positive definite of some appropriate
order. The name given to interpolants employing these basic
functions is surface splines. This is because they are a genuine
multivariate analogue of the well-loved natural splines in one
dimension.

It is of central importance to understand the behaviour of the error
between a function $f:\Omega \To \R$ and its interpolant as the set
$\nodes \subset \Omega$ becomes ``dense'' in $\Omega$. The measure
of density we employ is the {\it fill-distance} $h=\sup_{x \in
\Omega} \min_{a \in \nodes} \abs{x-a}$. One might hope that for some
suitable norm $\norm{\,\cdot\,}$ there is a constant $\gamma$,
independent of $f$ and $h$, such that
\begin{equation*}
\norm{f-Sf} = \order{(h^\gamma)},\qquad \mbox{as $h \To 0$}.
\end{equation*}
In the case of the Beppo-Levi spaces, there is a considerable
freedom of choice for the norm in which the error between $f$ and
$Sf$ is measured. The most widely quoted result concerns the norm
$\norm{\,\cdot\,}_{L_\infty(\Omega)}$, but for variety we prefer to
deal with the $L_p$-norm. To do this it is helpful to assume
$\Omega$ is a bounded domain, whose boundary is sufficiently smooth.
In this case
 there is a constant $C>0$, independent of
$f$ and $h$, such that for all $f \in BL^k(\Omega)$,
\begin{equation}\label{known error}
    \norm{f-Sf}_{L_p(\Omega)} \leq \biggl\{
        \begin{array}{ll}
            C h^{k-\frac{d}{2}+\frac{d}{p}}\abs{f}_{k,\Omega},& 2\leq p \leq
            \infty\\
            C h^{k} \abs{f}_{k,\Omega},& 1\leq p <2
        \end{array},\qquad \mbox{as $h \To 0$.}
\end{equation}

There has been considerable interest recently in the following very
natural question. What happens if the function $f$ does not possess
sufficient smoothness to lie in $BL^k(\Omega)$? It may well be that
$f$ lies in $BL^m(\Omega)$, where $2k>2m>d$. The condition $2m>d$
ensures that $f(a)$ exists for each $a\in\nodes$, and so $Sf$
certainly exists. However, $\abs{f}_{k,\Omega}$ is not defined. It
is simple to conjecture that the new error estimate should be
\begin{equation}\label{ourerror}
    \norm{f-Sf}_{L_p(\Omega)} \leq \biggl\{
        \begin{array}{ll}
            C h^{m-\frac{d}{2}+\frac{d}{p}}\abs{f}_{m,\Omega},& 2\leq p \leq
            \infty\\
            C h^{m} \abs{f}_{m,\Omega},& 1\leq p <2
        \end{array},\qquad \mbox{as $h \To 0$}.
\end{equation}
It is perhaps surprising to the uninitiated reader that this
estimate is not true even with the reasonable restrictions we have
placed on $k$ and $m$. We are going to describe a recent result
from~\citeasnoun{johnsonnew}. To do that, we recall the familiar
definition of a Sobolev space. Let $W^k_2(\Omega)$ denote the
$k$\textsuperscript{th} order Sobolev space, which consists of
functions all of whose derivatives up to and including order $k$ are
in $L_2(\Omega)$. It is a Banach space under the norm
\begin{equation*}
\norm{f}_{k,\Omega}=\Biggl( \sum_{i=0}^k \abs{f}_{i,\Omega}^2
\Biggr)^{1/2},\qquad \mbox{where $f\in W^k_2(\Omega)$}.
\end{equation*}
We have already tacitly alluded to the Sobolev embedding theorem
which states that when $\Omega$ is reasonably regular (for example,
when $\Omega$ possesses a Lipschitz continuous boundary) and
$k>d/2$, then the space $W^k_2(\Omega)$ can be embedded in
$C(\Omega)$~\citeaffixed[Theorem 5.4, p. 97]{Adams}{see}. Now
Johnson's result is as follows.
\begin{thm}[\citename{johnsonnew}]
Let $\Omega$ be the unit ball in $\R^d$ and assume $d/2<m<k$. For
every $h_0>0$, there exists an $f \in W^m_2(\R^d)$ and a sequence of
sets $\{\nodes_n\}_{n\in \N}$ with the following properties:
\begin{list}{}{}
\item (i) each set $\nodes_n$ consists of finitely many points
contained in $\Omega$; \item (ii) the fill-distance of each set
$\nodes_n$ is at most $h_0$; \item (iii) if $S^n_kf$ is the surface
spline interpolant to $f$ from $BL^k(\R^d)$ associated with
$\nodes_n$, for each $n\in \N$, then $\norm{S_k^nf}_{L_1(\Omega)}
\rightarrow \infty$ as $n\rightarrow \infty$.
\end{list}{}{}
\end{thm}
If the surface spline interpolation operator is unbounded, there is
of course no possibility of getting an error estimate of the kind we
conjectured. Johnson's proof uses point sets which have a special
feature. We define the separation distance of $\nodes_n$ as $q_n =
\min\{\abs{a-b}/2: a,b\in \nodes_n, a\neq b\}$. Let the
fill-distance of each $\nodes_n$ be $h_n$. In Johnson's proof, the
construction of $\nodes_n$ is such that $q_n/h_n\rightarrow 0$. We
make this remark, because Johnson's result in one dimension refers
to interpolation by natural splines, and in this setting the
connection between the separation distance and the unboundedness of
$S_k^n$ has been known for some time. What is also known in the
one-dimensional case is that if the separation distance is tied to
the fill-distance, then a result of the type we are seeking is true.
Theorem \ref{main} is the definitive result we obtain, and is the
formalisation of the conjectured bounds in
Equation~(\ref{ourerror}).

Subsequent to carrying out this work, we became aware of independent
work by~\citeasnoun{yoonnew}. In that paper, error bounds for the
case we consider here are also offered. Because of Yoon's technique
of proof, which is considerably different to our own, he obtains
error bounds for functions $f$ with the additional restriction that
$f$ lies in $W^k_\infty (\Omega)$, so the results here have wider
applicability. However, Yoon does consider the shifted surface
splines, whilst in this paper we have chosen to consider only
surface splines as an exemplar of what can be achieved. At the end
of Section~3 we offer some comments on the difference between our
approach and that of Yoon.

To close this section we introduce some notation that will be
employed throughout the paper. The support of a function $\phi:\R^d
\To \R$ is defined to be the closure of the set $\set{x \in \R^d:\
\phi(x) \neq 0}$, and is denoted by $\supp(\phi)$. The volume of a
bounded set $\Omega$ is the quantity $\int_\Omega \,dx$ and will be
denoted $\vol{\Omega}$. We make much use of the space $\Pi_{m-1}$,
so for brevity we fix $\ell$ as the dimension of this space.
Finally, when we write $\h{f}$ we mean the Fourier transform of $f$.
The context will clarify whether the Fourier transform is the
natural one on $L_1(\R^d)$:
$$
\widehat{f}(x) = \frac{1}{(2\pi)^{d/2}}\int_{\R^d} f(t) e^{-ixt}\,
dt,
$$
or one of its several extensions to $L_2(\R^d)$ or ${\cal S}'$.

%
%

\section{SOBOLEV EXTENSION THEORY}

In this section we intend to collect together a number of useful
results, chiefly about the sorts of extensions which can be carried
out on Sobolev spaces. We begin with the well-known result which can
be found in many of the standard texts. Of course, the precise
nature of the set $\Omega$ in the following theorem varies from book
to book, and we have not striven here for the utmost generality,
because that is not really a part of our agenda in this paper.

\begin{thm}[\citename{Adams}~\citeyear*{Adams}, Theorem 4.32, p. 91]\label{sobolev extension thm}
    Let $\Omega$ be an open, bounded subset of $\R^d$ satisfying the uniform
cone condition.
    For every $f \in W^{m}_2(\Omega)$ there
    is an $f^\Omega \in W^{m}_{2}(\R^d)$ satisfying
    $f^\Omega\!\!\mid_{\Omega} = f$. Moreover, there is a positive
    constant $K=K(\Omega)$ such that for all $f \in W^{m}_2(\Omega)$,
        \begin{equation*}
            \norm{f^\Omega}_{m,\R^d} \leq K \norm{f}_{m,\Omega}.
        \end{equation*}
\end{thm}

We remark that the extension $f^\Omega$ can be chosen to be
supported on any compact subset of $\R^d$ containing $\Omega$. To
see this, we construct $f^\Omega$ in accordance with
Theorem~\ref{sobolev extension thm}, then select $\eta \in
C^m_0(\R^d)$ such that $\eta(x)=1$ for $x \in \Omega$. Now, if we
consider the compactly supported function $f^\Omega_0  = \eta
f^{\Omega} \in W^m_2(\R^d)$, we have $f^\Omega_0 \!\!\mid_\Omega =
f$. An elementary application of the Leibniz formula gives
\begin{equation*}
    \norm{f^\Omega_0}_{m,\R^d} \leq C \norm{f}_{m,\Omega},\qquad
\mbox{where
    $C=C(\Omega,\eta)$}.
\end{equation*}

One of the nice features of the above extension is that the
behaviour of the constant $K(\Omega)$ can be understood for simple
choices of $\Omega$. The reason for this is of course the choice of
$\Omega$ and the way the seminorms defining the Sobolev norms behave
under dilations and translations of $\Omega$.

\begin{lemma}\label{cov lemma}
    Let $\Omega$ be a measurable subset of $\R^d$. Define the mapping
$\sigma:\R^d\rightarrow\R^d$ by
    $\sigma(x)=a+h(x-t)$, where $h>0$, and
    $a$, $t$, $x \in \R^d$. Then for all $f \in {W}^m_2(\sigma(\Omega))$,
        \begin{equation*}
            \abs{f \circ \sigma}_{m,\Omega} = h^{m-d/2}
\abs{f}_{m,\sigma(\Omega)}.
        \end{equation*}
\end{lemma}

\proof We have, for $\abs{\alpha}=m$,
\begin{equation*}
    (D^\alpha (f \circ \sigma))(x) = h^m (D^\alpha f )(\sigma(x)).
\end{equation*}
Thus,
\begin{equation*}
    \begin{split}
        \abs{f \circ \sigma}^2_{m,\Omega} &= \sum_{\abs{\alpha}=m}
c_\alpha
                                            \int_{\Omega} \abs{(D^\alpha
(f\circ \sigma))(x)}^2\
                                            dx\\
                                          &= h^{2m}
\sum_{\abs{\alpha}=m} c_\alpha
                                            \int_{\Omega} \abs{(D^\alpha
f )(\sigma(x))}^2\,
                                            dx.
    \end{split}
\end{equation*}
Now, using the change of variables $y=\sigma(x)$,
\begin{equation*}
    \abs{f \circ \sigma}^2_{m,\Omega} = h^{2m-d} \sum_{\abs{\alpha}=m}
c_\alpha
                                            \int_{\sigma(\Omega)}
\abs{(D^\alpha f )(y)}^2\,
                                            dy =
h^{2m-d}\abs{f}^2_{m,\sigma(\Omega)}. \qedhere
\end{equation*}

Unfortunately, the Sobolev extension refers to the Sobolev norm. We
want to work with a norm which is more convenient for our purposes.
This norm is in fact equivalent to the Sobolev norm, as we shall now
see.

\begin{lemma}\label{norm equivlance}
    Let $\Omega$ be an open subset of $\R^d$ having the cone property and a Lipschitz-continuous boundary. Let
    $\dotz{b_1}{b_\ell} \in \Omega$ be unisolvent with respect to
    $\Pi_{m-1}$. Define a norm on $W^m_2(\Omega)$ via
        \begin{equation*}
            \norm{f}_\Omega = \Biggl(
\abs{f}_{m,\Omega}^2+\sum_{i=1}^\ell
            \abs{f(b_i)}^2 \Biggr)^{1/2},\qquad f\in W^m_2(\Omega).
        \end{equation*}
    There are positive constants $K_1$ and
    $K_2$ such that for all $f \in W^m_2(\Omega)$,
        \begin{equation*}
            K_1 \norm{f}_{m,\Omega} \leq \norm{f}_{\Omega}  \leq K_2
            \norm{f}_{m,\Omega}.
        \end{equation*}
\end{lemma}

\proof The conditions imposed on $m$ and $\Omega$ ensure that
$W^m_2(\Omega)$ is continuously embedded in
$C(\Omega)$~\cite[Theorem 5.4, p. 97]{Adams}. So, given $x \in
\Omega$, there is a constant $C$ such that $\abs{f(x)} \leq C
\norm{f}_{m,\Omega}$ for all $ f \in W^m_2(\Omega)$. Thus, there are
constants $\dotz{C_1}{C_\ell}$ such that
    \begin{equation}
        \norm{f}_{\Omega}^2 \leq \abs{f}_{m,\Omega}^2 + \sum_{i=1}^\ell
        C_i \norm{f}_{m,\Omega}^2 \leq \Big(1+\sum_{i=1}^\ell C_i\Big)
        \norm{f}_{m,\Omega}^2.
\label{bound}
    \end{equation}
On the other hand, suppose there is no positive number $K$ with
$\norm{f}_{m,\Omega} \leq K \norm{f}_{\Omega}$ for all $f \in
W^m_2(\Omega)$. Then there is a sequence $\set{f_j}$ in
$W^m_2(\Omega)$ with
    \begin{equation*}
        \norm{f_j}_{m,\Omega} =1\qquad \mbox{and}\qquad
        \norm{f_j}_{\Omega} \leq \frac{1}{j},\qquad \mbox{for $j=1,2,\dotsc$.}
    \end{equation*}
The Rellich selection theorem~\cite[Theorem 1.9, p. 32]{Braess}
states that $W^m_2(\Omega)$ is compactly embedded in
$W^{m-1}_2(\Omega)$. Therefore, as $\set{f_j}$ is bounded in
$W^m_2(\Omega)$, this sequence must contain a convergent subsequence
in $W^{m-1}_2(\Omega)$. With no loss of generality we shall assume
$\set{f_j}$ itself converges in $W^{m-1}_2(\Omega)$. Thus $\{f_j\}$
is a Cauchy sequence in $W^{m-1}_2(\Omega)$. Next, as
$\norm{f_j}_{\Omega} \To 0$ it follows that $\abs{f_j}_{m,\Omega}
\To 0$. Moreover,
\begin{equation*}
    \begin{split}
        \norm{f_j - f_k}_{m,\Omega}^2 &= \norm{f_j - f_k}_{m-1,\Omega}^2
+
                                        \abs{f_j -
                                        f_k}_{m,\Omega}^2\\
                                      &\leq \norm{f_j -
f_k}_{m-1,\Omega}^2 +
                                        2\abs{f_j}_{m,\Omega}^2 +
                                        2\abs{f_k}_{m,\Omega}^2.
        \end{split}
\end{equation*}
Since $\{f_j\}$ is a Cauchy sequence in $W^{m-1}_2(\Omega)$, and
$\abs{f_j}_{m,\Omega} \rightarrow 0$, it follows that $\set{f_j}$ is
a Cauchy sequence in $W^{m}_2(\Omega)$.  Since $W^m_2(\Omega)$ is
complete with respect to $\norm{\,\cdot\,}_{m,\Omega}$, this
sequence converges to a limit $f \in W^m_2(\Omega)$. By
Equation~(\ref{bound}),
$$
\norm{f-f_j}_\Omega^2 \leq \Big(1+\sum_{i=1}^\ell C_i\Big)
\norm{f-f_j}^2_{m,\Omega},
$$
and hence $\norm{f-f_j}_\Omega \rightarrow 0$ as $j\rightarrow
\infty$. Since $\norm{f_j}_\Omega \rightarrow 0$, it follows that
$f=0$. Because $\norm{f_j}_{m,\Omega} =1$, $j=1,2,\ldots$, it
follows that $\norm{f}_{m,\Omega}=1$. This contradiction establishes
the result. \qed

We are almost ready to state the key result which we will employ in
our later proofs about error estimates. Before we do this, let us
make a simple observation. Look at the unisolvent points
$\dotz{b_1}{b_\ell}$ in the statement of the previous Lemma. Since
$W^m_2(\Omega)$ can be embedded in $C(\Omega)$, it makes sense to
talk about the interpolation projection $P:W^m_2(\Omega) \rightarrow
\Pi_{m-1}$ based on these points. Furthermore, under certain nice
conditions (for example $\Omega$ being a bounded domain), $P$ is the
orthogonal projection of $W^m_2(\Omega)$ onto $\Pi_{m-1}$.

\begin{lemma}\label{poly extention}
    Let $B$ be any ball of radius $h$ and center $a \in \R^d$, and let
    $f \in W^{m}_{2}(B)$. Whenever $\dotz{b_1}{b_\ell} \in \R^d$ are unisolvent
    with respect to $\Pi_{m-1}$ let $P_b:C(\R^d) \To \Pi_{m-1}$ be the
    Lagrange interpolation operator on $\dotz{b_1}{b_\ell}$. Then
    there exists $c=(\dotz{c_1}{c_\ell}) \in B^\ell$ and
    $g \in W^{m}_{2}(\R^d)$ such that
        \begin{enumerate}
            \item $g(x)=(f-P_c f)(x)$ for all $x \in B$;

            \item $g(x)=0$ for all $\abs{x-a}>2h$;

            \item there exists a $C>0$, independent of $f$ and $B$, such
            that $\abs{g}_{m,\R^d} \leq C \abs{f}_{m,B}$.
        \end{enumerate}
    Furthermore, $\dotz{c_1}{c_\ell}$ can be arranged so that $c_1=a$.
\end{lemma}

\proof Let $B_1$ be the unit ball in $\R^d$ and let $B_2=2B_1$. Let
$\dotz{b_1}{b_\ell} \in B_1$ be unisolvent with respect to
$\Pi_{m-1}$. Define $\sigma(x)=h^{-1}(x-a)$ for all $x \in \R^d$.
Set $c_i = \sigma^{-1}(b_i)$ for $i=\dotz{1}{\ell}$ so that
$\dotz{c_1}{c_\ell} \in B$ are unisolvent with respect to
$\Pi_{m-1}$. Take $f \in W^m_2(B)$. Then $(f-P_c f) \circ
\sigma^{-1} \in W^m_2(B_1)$. Set $F=(f-P_c f) \circ \sigma^{-1}$.
Let $F^{B_1}$ be constructed as an extension to $F$ on $B_1$. By
Theorem~\ref{sobolev extension thm} and the remarks following it, we
can assume $F^{B_1}$ is supported on $B_2$. Define $g=F^{B_1} \circ
\sigma \in W^m_2(\R^d)$. Let $x \in B$. Since $\sigma(B)=B_1$ there
is a $y \in B_1$ such that $x=\sigma^{-1}(y)$. Then,
\begin{equation*}
    g(x)=(F^{B_1} \circ \sigma)(x) = F^{B_1}(y) = ((f-P_c f) \circ
    \sigma^{-1}) (y) = (f-P_c f)(x).
\end{equation*}
Also, for $x \in \R^d$ with $\abs{x-a}>2h$, we have $\abs{\sigma(x)}
> 2$. Since $F^{B_1}$ is supported on $B_2$, $g(x)=0$ for $\abs{x-a}>2h$. Hence,
$g$ satisfies properties \textit{1} and \textit{2}. By
Theorem~\ref{sobolev extension thm} there is a $K_1$, independent of
$f$ and $B$, such that
\begin{equation*}
    \norm{F^{B_1}}_{m,B_2} = \norm{F^{B_1}}_{m,\R^d} \leq K_1
    \norm{F}_{m,B_1}.
\end{equation*}
We have seen in Lemma~\ref{norm equivlance} that if we endow
$W^m_2(B_1)$ and $W^m_2(B_2)$ with the norms
\begin{equation*}
    \norm{v}_{B_i} = \Biggl( \abs{v}_{m,B_i}^2 + \sum_{i=1}^\ell
    \abs{v(b_i)}^2 \Biggr)^{1/2},
\qquad i=1,2,
\end{equation*}
then $\norm{\,\cdot\,}_{B_i}$ and $\norm{\,\cdot\,}_{m,B_i}$ are
equivalent for $i=1,2$.  Thus, there are constants $K_2$ and $K_3$,
independent of $f$ and $B$, such that
\begin{equation*}
    \norm{F^{B_1}}_{B_2} \leq K_2 \norm{F^{B_1}}_{m,B_2} \leq K_1 K_2
    \norm{F}_{m,B_1} \leq K_1 K_2 K_3 \norm{F}_{B_1}.
\end{equation*}
Set $C=K_1 K_2 K_3$. Since $F^{B_1}(b_i)=F(b_i)=(f-P_c
f)(\sigma^{-1}(b_i))=(f-P_c f)(c_i)=0$ for $i=\dotz{1}{\ell}$, it
follows that $\abs{F^{B_1}}_{m,B_2} \leq C \abs{F}_{m,B_1}$. Thus,
$\abs{g \circ \sigma^{-1}}_{m,\R^d} \leq C \abs{(f-P_c f) \circ
\sigma^{-1}}_{m,B_1}$. Now, Lemma~\ref{cov lemma} can be employed
twice to give
\begin{equation*}
    \abs{g}_{m,\R^d} = h^{d/2 -m}\abs{g \circ \sigma^{-1}}_{m,\R^d} \leq C h^{d/2
    -m} \abs{(f-P_c f) \circ \sigma^{-1}}_{m,B_1} = C \abs{f-P_c f}_{m,B}.
\end{equation*}
Finally, we observe that $\abs{f-P_c f}_{m,B} = \abs{f}_{m,B}$ to
complete the first part of the proof. The remaining part follows by
selecting $b_1=0$ and choosing $\dotz{b_2}{b_\ell}$ accordingly in
the above construction. \qed

\begin{lemma}[\citename{Duchon2}~\citeyear*{Duchon2}]\label{light wayne ext}
Let $\Omega$ be an open, bounded, connected subset of $\R^d$ having
the cone property and a Lipschitz-continuous boundary. Let $f \in
W^m_2(\Omega)$. Then there exists a unique element $f^\Omega \in
BL^m(\R^d)$ such that $f^\Omega\!\!\mid_{\Omega} = f$, and amongst
all elements of $BL^m(\R^d)$ satisfying this condition,
$\abs{f^\Omega}_{m,\R^d}$ is minimal. Furthermore, there exists a
constant $K=K(\Omega)$ such that, for all $f \in W^m_2(\Omega)$,
\begin{equation*}
\abs{f^\Omega}_{m,\R^d} \leq K \abs{f}_{m,\Omega}.
\end{equation*}
\end{lemma}

%
%

\section{ERROR ESTIMATES}

We arrive now at our main section, in which we derive the required
error estimates. Our strategy is simple. We begin with a function
$f$ in $BL^{m}(\R^d)$. We want to estimate $\norm{f- S_kf}$ for some
suitable norm $\norm{\,\cdot\,}$, where $S_k$ is the minimal norm
interpolation operator from $BL^{k}(\R^d)$, and $k> m$. We suppose
that we already have an error bound using the norm
$\norm{\,\cdot\,}$ for all functions  $g\in BL^{k}(\R^d)$. Our proof
now proceeds as follows. Firstly, we adjust $f$ in a somewhat
delicate manner, obtaining a function $F$, still in $BL^{m}(\R^d)$,
and with seminorm in $BL^{m}(\R^d)$ not too far away from that of
$f$. We then smooth $F$ by convolving it with a function $\phi \in
C^\infty_0(\R^d)$. The key feature of the adjustment of $f$ to $F$
is that $(\phi*F)(a) = f(a)$ for every point $a$ in our set of
interpolation points. It then follows that $F\in BL^{k}(\R^d)$. We
then use the usual error estimate in $BL^{k}(\R^d)$. A standard
procedure (Lemma~\ref{seminorm convolution bound}) then takes us
back to an error estimate in $BL^{m}(\R^d)$.

\begin{lemma}\label{seminorm convolution bound}
    Let $m \leq k$ and let $\phi \in C^{\infty}_0(\R^{d})$. For each $h>0$, let
    $\phi_h(x) = h^{-d} \phi(x/h)$ for $x \in \R^d$. Then there exists
    a constant $C>0$, independent of $h$, such that for all $f \in
    {BL}^{m}(\R^d)$,
        \begin{equation*}
            \abs{\phi_h * f}_{k,\R^d} \leq C h^{m-k}
\abs{f}_{m,\R^d}.
        \end{equation*}
    Furthermore, we have $\abs{\phi_h *
    f}_{k,\R^d} = o(h^{m-k})$ as $h \To 0$.
\end{lemma}

\proof The chain rule for differentiation gives $(D^\gamma
\phi_h)(x) = h^{-(d+\abs{\gamma})} (D^\gamma \phi) (x/h)$ for all $x
\in \R^d$, and $\gamma \in \Z^d_+$. Thus, for $\beta \in \Z^d_+$
with $\abs{\beta}=m$ we have
\begin{align}
        \int_{\R^d} \abs{(D^\gamma \phi_h * D^\beta f)(x)}^2\ dx
                        &= \int_{\R^d} \biggabs{ \int_{\R^d} (D^\gamma
\phi_h) (x-y)
                           (D^\beta f)(y)\ dy}^2\ dx \nonumber\\
                        &= h^{-2(d+\abs{\gamma})} \int_{\R^d} \biggabs{
\int_{\R^d}
                           (D^\gamma \phi) \Bigl(\frac{x-y}{h}\Bigr)
(D^\beta f)(y)\, dy}^2\
                            dx \nonumber\\
                        &= h^{-2\abs{\gamma}} \int_{\R^d} \biggabs{
\int_{\R^d} (D^\gamma \phi )(t)
                           (D^\beta f)(x-h t)\ dt}^2\, dx \nonumber\\
                        &= h^{-2\abs{\gamma}} \int_{\R^d} \biggabs{
\int_{K} (D^\gamma \phi )(t)
                           (D^\beta f)(x-h t)\ dt}^2\, dx,\label{little
oh}
\end{align}
where $K=\supp{(\phi)}$. An application of the Cauchy-Schwartz
inequality gives
\begin{equation*}
    \int_{\R^d} \abs{(D^\gamma \phi_h * D^\beta f)(x)}^2\ dx \leq
    h^{-2\abs{\gamma}} \int_{\R^d} \Biggl( \int_{K} \abs{(D^\gamma
    \phi)(t)}^2\ dt \Biggr)\Biggl(\int_{K} \abs{(D^\beta f)(x-h t)}^2\
    dt\Biggr)\ dx,
\end{equation*}
and so,
\begin{equation}\label{post cauchy-schwartz}
    \int_{\R^d} \abs{(D^\gamma \phi_h * D^\beta f)(x)}^2\ dx \leq
    h^{-2\abs{\gamma}} \int_{\R^d} \abs{(D^\gamma \phi)(t)}^2\ dt
    \int_{\R^d} \int_{K} \abs{(D^\beta f)(x-h t)}^2\ dtdx.
\end{equation}
The Parseval formula together with the relation $(D^\alpha (\phi_h
* f))\ehat{} = (i\,\cdot\,)^\alpha (\phi_h
    * f)\ehat{}$ provide us with the
equality
\begin{align}
            \sum_{\abs{\alpha}=k} c_\alpha \int_{\R^d} \abs{(D^\alpha
(\phi_h * f))(x)}^2\ dx
                    &= \sum_{\abs{\alpha}=k} c_\alpha \int_{\R^d}
\abs{(ix)^\alpha (\phi_h * f)\ehat(x)}^2\
                    dx\nonumber\\
                    &= \int_{\R^d} \sum_{\abs{\alpha}=k} c_\alpha
x^{2\alpha} \abs{(\phi_h * f)\ehat(x)}^2\
                    dx.\label{i x alpha}
\end{align}
Now, when Equation \eqref{i x alpha} is used in conjunction with the
relation
\begin{equation*}
    \sum_{\abs{\alpha}=k} c_\alpha x^{2\alpha} = \abs{x}^{2k}=
    \abs{x}^{2(m+k-m)}= \sum_{\abs{\beta}=m} c_\beta x^{2\beta}
    \sum_{\abs{\gamma}=k-m} c_\gamma x^{2\gamma},
\end{equation*}
we obtain
\begin{equation*}
    \begin{split}
        \sum_{\abs{\alpha}=k} c_\alpha \int_{\R^d} \abs{(D^\alpha
(\phi_h
        * f))(x)}^2\ dx             &= \int_{\R^d}
                           \sum_{\abs{\beta}=m} c_\beta x^{2\beta}
                           \sum_{\abs{\gamma}=k-m} c_\gamma x^{2\gamma}
                           \abs{(\phi_h * f)\ehat(x)}^2\ dx\\
                         &= \sum_{\abs{\beta}=m} c_\beta \int_{\R^d}
                            \sum_{\abs{\gamma}=k-m} c_\gamma x^{2\gamma}
                            \abs{(ix)^\beta (\phi_h * f)\ehat(x)}^2\ dx
\\
                         &= \sum_{\abs{\beta}=m} c_\beta \int_{\R^d}
                            \sum_{\abs{\gamma}=k-m} c_\gamma x^{2\gamma}
                            \abs{(D^\beta(\phi_h * f))\ehat(x)}^2\
dx \\
                         &= \sum_{\abs{\beta}=m} c_\beta
\sum_{\abs{\gamma}=k-m} c_\gamma
                            \int_{\R^d}
                            \abs{(ix)^{\gamma} (D^\beta(\phi_h *
f))\ehat(x)}^2\ dx \\
                         &= \sum_{\abs{\beta}=m} c_\beta
\sum_{\abs{\gamma}=k-m} c_\gamma
                            \int_{\R^d}
                            \abs{(D^\gamma (D^\beta(\phi_h *
f)))\ehat(x)}^2\ dx \\
                         &= \sum_{\abs{\beta}=m} c_\beta
\sum_{\abs{\gamma}=k-m} c_\gamma
                            \int_{\R^d}
                            \abs{(D^\gamma (D^\beta(\phi_h *
f)))(x)}^2\ dx .
    \end{split}
\end{equation*}
Since the operation of differentiation commutes with convolution, we
have that
\begin{equation}\label{pre cauchy-schwarz}
    \sum_{\abs{\alpha}=k} c_\alpha \int_{\R^d} \abs{(D^\alpha (\phi_h
    * f))(x)}^2\ dx = \sum_{\abs{\beta}=m} c_\beta
    \sum_{\abs{\gamma}=k-m} c_\gamma \int_{\R^d} \abs{(D^\gamma \phi_h
    * D^\beta f )(x)}^2\ dx.
\end{equation}
Combining Equation \eqref{post cauchy-schwartz} with Equation
\eqref{pre cauchy-schwarz} we deduce that
\begin{align*}
    \sum_{\abs{\alpha}=k} c_\alpha \int_{\R^d} &\abs{(D^\alpha (\phi_h
    * f))(x)}^2\ dx
 \\
    &\leq  \sum_{\abs{\beta}=m} c_\beta \sum_{\abs{\gamma}=k-m} c_\gamma
h^{-2\abs{\gamma}} \int_{\R^d} \abs{(D^\gamma \phi)(t)}^2\ dt
\int_{\R^d} \int_{K} \abs{(D^\beta f)(x-h t)}^2\ dtdx\\
    &=  h^{2(m-k)} \abs{\phi}_{k-m,\R^d}^2 \sum_{\abs{\beta}=m} c_\beta
\int_{\R^d} \int_{K} \abs{(D^\beta f)(x-h t)}^2\ dtdx.
\end{align*}
Fubini's theorem permits us to change the order of integration in
the previous inequality. Thus,
\begin{equation*}
    \sum_{\abs{\alpha}=k} c_\alpha \int_{\R^d} \abs{(D^\alpha (\phi_h
    * f))(x)}^2\, dx \leq h^{2(m-k)} \abs{\phi}_{k-m,\R^d}^2
    \sum_{\abs{\beta}=m} c_\beta \int_{K} \int_{\R^d} \abs{(D^\beta
    f)(x-h t)}^2\, dxdt.
\end{equation*}
Finally, a change of variables in the inner integral above yields
\begin{equation*}
    \sum_{\abs{\alpha}=k} c_\alpha \int_{\R^d} \abs{(D^\alpha (\phi_h
    * f))(x)}^2\, dx \leq h^{2(m-k)} \abs{\phi}_{k-m,\R^d}^2
    \sum_{\abs{\beta}=m} c_\beta \int_{K} \int_{\R^d} \abs{(D^\beta
    f)(z)}^2\, dzdt.
\end{equation*}
Setting $C=\abs{\phi}_{k-m,\R^d} \sqrt{\vol{K}}$ we conclude that
$\abs{\phi_h * f}_{k,\R^d} \leq C h^{m-k} \abs{f}_{m,\R^d}$ as
required. To deal with the remaining statement of the lemma, we
observe that for $\gamma \neq 0$ we have
\begin{equation*}
 \int_{K} (D^\gamma \phi)(t)\, dt= \int_{\R^d} (D^\gamma \phi)(t)\, dt = (\widehat{D^\gamma \phi})(0) =
((i\,\cdot\,)^\gamma \widehat{\phi})(0) = 0.
\end{equation*}
Then it follows from Equation~\eqref{little oh} that for
$\abs{\beta}=m$,
\begin{equation*}
     \int_{\R^d} \abs{(D^\gamma \phi_h * D^\beta f)(x)}^2\ dx =
            h^{-2\abs{\gamma}} \int_{\R^d} \biggabs{ \int_{K} (D^\gamma
\phi
            )(t)(  (D^\beta f)(x-h t) - (D^\beta f)(x)) \ dt}^2\ dx.
\end{equation*}
Now, if we continue in precisely the same manner as before, we
obtain
\begin{multline*}
    \sum_{\abs{\alpha}=k} c_\alpha \int_{\R^d} \abs{(D^\alpha (\phi_h
    * f))(x)}^2\ dx \\\leq h^{2(m-k)} \abs{\phi}_{k-m,\R^d}^2
    \sum_{\abs{\beta}=m} c_\beta \int_{K} \int_{\R^d} \abs{(D^\beta
    f)(x-h t) - (D^\beta f)(x)}^2\ dxdt.
\end{multline*}
Since $D^\beta f \in L^2(\R^d)$ for each $\beta \in \Z^d_+$ with
$\abs{\beta}=m$, it follows that for almost all $t,x \in \R^d$,
\begin{equation*}
    \abs{(D^\beta f)(x -ht) - (D^\beta f)(x)} \To 0,\qquad \mbox{as $h \To 0$.}
\end{equation*}
Furthermore, setting
\begin{equation*}
g(x,t) = 2\abs{(D^\beta f)(x-ht)}^2 + 2 \abs{(D^\beta f)(x)}^2,
\qquad \mbox{for almost all $x,t\in\R^d$,}
\end{equation*}
we see that
\begin{equation*}
\abs{(D^\beta f )(x -ht) - (D^\beta f)(x)}^2 \leq g(x,t),
\end{equation*}
for almost all $x,t\in\R^d$ and each $h>0$. It follows by
calculations similar to those used above that
\begin{equation*}
\int_{K} \int_{\R^d} g(x,t)\, dxdt =
4\mbox{vol}(K)\int_{\R^d}\abs{(D^\beta f)(x)}^2 \, dx <\infty.
\end{equation*}
Applying Lebesgue's dominated convergence theorem, we obtain
\begin{equation*}
\int_{K} \int_{\R^d} \abs{(D^\beta f)(x-h t) - (D^\beta f)(x)}^2\
dxdt \To 0,\qquad \mbox{as $h \To 0$.}
\end{equation*}
Hence, for $m\leq k$, $\abs{\phi_h *
    f}_{k,\R^d} = o(h^{m-k})$ as $h \To 0$. \qed

\begin{lemma}\label{poly repro}
    Suppose $\phi \in C^{\infty}_0(\R^{d})$ is supported on the unit ball and
    satisfies
    \begin{equation*}
        \int_{\R^d} \phi(x)\ dx = 1\qquad \mbox{and}\qquad \int_{\R^d}
        \phi(x) x^\alpha\ dx =0,\qquad \mbox{for all}\  0< \abs{\alpha}
\leq
        {m-1}.
    \end{equation*}
    For each $\eps>0$ and $x \in \R^d$, let $\phi_\eps(x) =
\eps^{-d}
    \phi(x/\eps)$. Let $B$ be any ball of radius $h$ and center $a \in
    \R^d$. For a fixed $p \in \Pi_{m-1}$ let $f$ be a mapping from
    $\R^d$ to $\R$ such that $f(x) = p(x)$ for all $x \in B$. Then
    $(\phi_\eps * f) (a) = p(a)$ for all $\eps \leq h$.
\end{lemma}

\proof Let $B_1$ denote the unit ball in $\R^d$. We begin by
employing a change of variables to deduce
\begin{equation*}
    \begin{split}
        (\phi_\eps * f) (a) &= \int_{\R^d} \phi_\eps(a-y) f(y)\ dy\\
                             &= \eps^{-d} \int_{\R^d}
                             \phi\Bigl(\frac{a-y}{\eps}\Bigr)
                             f(y)\ dy\\
                             &= \int_{\R^d} \phi(x)
                             f(a-x\eps)\ dx\\
                             &= \int_{B_1} \phi(x)
                             f(a-x\eps)\ dx.
    \end{split}
\end{equation*}
 Then, for $x \in B_1$,
$\abs{(a-x\eps)-a} \leq \eps \leq h$. Thus, $f(a-x\eps)=p(a-x\eps)$
for all $x\in B_1$. Moreover, there are numbers $b_\alpha$ such that
$p(a-x\eps)=p(a)+\sum_{0<\abs{\alpha}\leq {m-1}} b_\alpha x^\alpha$.
Hence,
\begin{equation*}
    \begin{split}
        (\phi_\eps * f) (a) &= \int_{B_1} \phi(x)
                             p(a-x\eps)\ dx\\
                             &= \int_{\R^d} \phi(x)
                             \biggl(p(a)+\sum_{0<\alpha\leq {m-1}}b_\alpha
x^\alpha \biggr)\ dx\\
                             &= p(a). \qedhere
    \end{split}
\end{equation*}
\begin{defin}
    Let $\Omega$ be an open, bounded subset of $\R^d$. Let $\nodes$ be a
    set of points in $\Omega$. The quantity $\sup_{x \in
    \Omega} \inf_{a \in \nodes} \abs{x-a} = h$ is called the fill-distance of
$\nodes$ in $\Omega$. The separation of $\nodes$ is
    given by the quantity
    \begin{equation*}
        q=\min_{\substack{a,b \in \nodes \\ a\neq b} }
        \frac{\abs{a-b}}{2}.
    \end{equation*}
The quantity $h/q$ will be called the mesh-ratio of $\nodes$.
\end{defin}
\begin{thm}\label{general result}
    Let $\nodes$ be a finite subset of $\R^d$ of separation $q>0$ and let $d<2m\leq 2k$. Then
    for all $f \in BL^m(\R^d)$ there exists an $F \in BL^k(\R^d)$
    such that
        \begin{enumerate}
            \item $F(a)=f(a)$ for all $a \in \nodes$;
            \item there exists a $C>0$, independent of $f$ and $q$, such that $\abs{F}_{m,\R^d}\leq C \abs{f}_{m,\R^d}$
             and $\abs{F}_{k,\R^d} \leq C q^{m-k}\abs{f}_{m,\R^d}$.
        \end{enumerate}
\end{thm}

\proof Take $f \in BL^m(\R^d)$. For each $a \in \nodes$ let $B_a
\subset \R^d$ denote the ball of radius $\delta = q/4$ centered at
$a$. For each $B_a$ let $g_a$ be constructed in accordance with
Lemma~\ref{poly extention}.  That is, for each $a \in \nodes$ take
$c'=(\dotz{c_2}{c_\ell}) \in B_a^{\ell-1}$ and $g_a \in W^m_2(\R^d)$
such that
\begin{enumerate}
    \item $a,\dotz{c_2}{c_\ell}$ are unisolvent with respect to
    $\Pi_{m-1}$;

    \item $g_a(x)=(f-P_{(a,c')} f)(x)$ for all $x \in B_a$;

    \item $P_{(a,c')} f \in \Pi_{m-1}$ and $(P_{(a,c')} f)(a) = f(a)$;

    \item $g_a(x)=0$ for all $\abs{x-a}>2\delta$;

    \item there exists a $C_1>0$, independent of $f$ and $B_a$, such
    that $\abs{g_a}_{m,\R^d} \leq C_1 \abs{f}_{m,B_a}$.
\end{enumerate}
Note that if $a\neq b$, then $\supp(g_a)$ does not intersect
$\supp(g_b)$, because if $x\in \supp(g_a)$ then
\begin{equation*}
    \abs{x - b} > \abs{b-a} - \abs{x-a} \geq 2q - 2 \delta = 6 \delta.
\end{equation*}
Using the observation above regarding the supports of the $g_a$'s it
follows that
\begin{align}
\biggabs{\sum_{a \in \nodes} g_a}_{m,\R^d}^2
        &=  \sum_{\abs{\alpha}=m} c_\alpha \int_{\R^d} \biggabs{\sum_{a \in
            \nodes} (D^\alpha g_a)(x)}^2\ dx \nonumber \\
        &=  \sum_{\abs{\alpha}=m} c_\alpha \sum_{b \in \nodes} \int_{\supp(g_b)} \biggabs{\sum_{a \in
            \nodes} (D^\alpha g_a)(x)}^2\ dx \nonumber \\
        &=  \sum_{\abs{\alpha}=m} c_\alpha \sum_{b \in \nodes} \int_{\supp(g_b)} \abs{(D^\alpha g_b)(x)}^2\ dx \nonumber \\
        &=  \sum_{a \in \nodes} \abs{g_a}_{m,\R^d}^2. \nonumber
\end{align}
Applying Condition 5 to the above equality we have
\begin{equation*}
    \biggabs{\sum_{a \in \nodes} g_a}_{m,\R^d}^2 \leq C_1^2 \sum_{a \in \nodes}
    \abs{f}_{m,B_a}^2\leq C_1^2 \abs{f}_{m,\R^d}^2.
\end{equation*}
Now set $H= f-\sum_{a\in \nodes} g_a$. It then follows from
Condition 2 above that $H(x) = (P_{(a,c')} f)(x)$ for all $x\in
B_a$, and from Condition 3 that $H(a) = f(a)$ for all $a\in \nodes$.
Let $\phi \in C^{\infty}_0(\R^{d})$ be supported on the unit ball
and enjoy the properties
\begin{equation*}
    \int_{\R^d} \phi(x)\ dx = 1\qquad \mbox{and}\qquad \int_{\R^d}
    \phi(x) x^\alpha\ dx =0,\qquad \mbox{for all}\ 0< \abs{\alpha} \leq
    {m-1}.
\end{equation*}
Now set $F=\phi_\delta * H$. Using Lemma~\ref{seminorm convolution
bound}, there is a constant $C_2>0$, independent of $q$ and $f$,
such that
\begin{align}
        \abs{F}_{k,\R^d}^2
            &\leq   C_2 \delta^{2(m-k)} \biggabs{f - \sum_{a \in \nodes} g_a}_{m,\R^d}^2 \nonumber\\
            &\leq   2 C_2 \delta^{2(m-k)} \biggl( \abs{f}_{m,\R^d}^2+ \biggabs{\sum_{a \in \nodes}
            g_a}_{m,\R^d}^2 \biggr)\nonumber\\
            &\leq   2 C_2 (1+C_1^2) \delta^{2(m-k)} \abs{f}_{m,\R^d}^2.
            \nonumber
\end{align}
Similarly, there is a constant $C_3>0$, independent of $q$ and $f$,
such that
\begin{align}
        \abs{F}_{m,\R^d}^2
            &\leq   C_3 \biggabs{f - \sum_{a \in \nodes} g_a}_{m,\R^d}^2 \nonumber\\
            &\leq   2 C_3 (1+C_1^2)\abs{f}_{m,\R^d}^2.
            \nonumber
\end{align} Thus
$\abs{F}_{k,\R^d} \leq C q^{m-k} \abs{f}_{m,\R^d}$ and
$\abs{F}_{m,\R^d} \leq  C \abs{f}_{m,\R^d}$ for some appropriate
constant $C>0$. Since $F=\phi_\delta
* H$ and $H|_{B_a}\in \Pi_{m-1}$ for each $a\in \nodes$, it
follows from Lemma~\ref{poly repro} that $F(a)=H(a) =f(a)$ for all
$a \in \nodes$.\qed

\begin{thm}
\label{main}
    Let $\Omega$ be an open, bounded, connected subset of $\R^d$ satisfying the
cone property and having a Lipschitz-continuous
    boundary.  Suppose also $d<2m\leq 2k$. For each
    $h>0$, let $\nodes_h$ be a finite,
    $\Pi_{k-1}$--unisolvent subset of $\Omega$ with fill-distance $h$. Assume also that there is a
quantity $\rho>0$ such that the mesh-ratio of each $\nodes_h$ is
bounded by $\rho$ for all $h>0$. For each mapping
    $f:\nodes_h\rightarrow \R$, let $S_k^h f$ be the
    minimal norm interpolant to $f$ on $\nodes_h$ from $BL^k (\R^d)$.
    Then there exists a constant $C>0$, independent of
    $h$, such that for all $f \in BL^m(\Omega)$,
    \begin{equation*}
        \norm{f-S^h_k f}_{L_p(\Omega)} \leq \biggl\{
            \begin{array}{ll}
                C h^{m-\frac{d}{2}+\frac{d}{p}}\abs{f}_{m,\Omega},& 2\leq p \leq
                \infty\\
                C h^{m} \abs{f}_{m,\Omega},& 1\leq p <2
            \end{array},\qquad \mbox{as $h \To 0$}.
    \end{equation*}
\end{thm}

\proof Take $f\in BL^m(\Omega)$. By~\citeasnoun{Duchon}, $f\in
W^m_2(\Omega)$. We define $f^\Omega$ in accordance with
Lemma~\ref{light wayne ext}. For most of this proof we wish to work
with $f^\Omega$ and not $f$, so for convenience we shall write $f$
instead of $f^\Omega$. Construct $F$ in accordance with
Theorem~\ref{general result} and set $G=f-F$. Then $F(a)=f(a)$ and
$G(a)=0$ for all $a \in \nodes_h$. Furthermore, there is a constant
$C_1>0$, independent of $f$ and $h$, such that
\begin{align}
&\abs{F}_{k,\R^d} \leq C_1 \left(\frac{h}{\rho}\right)^{m-k} \abs{f}_{m,\R^d},\label{F bound}\\
&\abs{G}_{m,\R^d} \leq \abs{f}_{m,\R^d}+\abs{F}_{m,\R^d} \leq
(1+C_1)\abs{f}_{m,\R^d}.\label{G bound}
\end{align}
Thus $S_k^h f = S_k^h F$ and $S_m^h G = 0$, where we have adopted
the obvious notation for $S_m^h$. Hence,
\begin{equation*}
    \norm{f-S_k^h f}_{L_p(\Omega)} = \norm{f-S_k^hF}_{L_p(\Omega)} =
    \norm{F+G-S_k^hF}_{L_p(\Omega)} \leq \norm{F-S_k^hF}_{L_p(\Omega)} +
    \norm{G-S_m^hG}_{L_p(\Omega)}.
\end{equation*}
Now, employing \possessivecite{Duchon2} error estimates for surface
splines \eqref{known error}, there are positive constants $C_2>0$
and $C_3>0$, independent of $h$ and $f$, such that
\begin{equation*}
    \norm{f-S_k^h f}_{L_p(\Omega)} \leq C_2 h^{\beta(k)}
\abs{F}_{k,\Omega}
    + C_3 h^{\beta(m)} \abs{G}_{m,\Omega},\qquad \mbox{as $h\rightarrow 0$,}
\end{equation*}
where we have defined
\begin{equation*}
        \beta(j)=\biggl\{
            \begin{array}{ll}
                j-\frac{d}{2}+\frac{d}{p},& 2\leq p \leq
                \infty\\
                j,& 1\leq p <2
            \end{array}.
    \end{equation*}
Finally, using the bounds in Equations \eqref{F bound} and \eqref{G
bound} we have
\begin{equation*}
    \norm{f-S_k^h f}_{L_p(\Omega)} \leq C_4 h^{\beta(m)}
    \abs{f}_{m,\R^d},\qquad \mbox{as $h\rightarrow 0$,}
\end{equation*}
for some appropriate $C_4>0$. To complete the proof we remind
ourselves that we have substituted $f^\Omega$ with $f$, and so an
application of Lemma~\ref{light wayne ext} shows that we can find
$C_5>0$ such that
\begin{equation*}
\norm{f-S_k^h f}_{L_p(\Omega)} \leq C_4 h^{\beta(m)}
    \abs{f^\Omega}_{m,\R^d} \leq C_4 C_5 h^{\beta(m)}
    \abs{f}_{m,\Omega},\qquad \mbox{as $h\rightarrow 0$.} \qedhere
\end{equation*}

We conclude this section with a brief commentary on the approach
of~\citeasnoun{yoonnew}. It is hardly surprising that Yoon's
technique also utilises a smoothing via convolution with a smooth
kernel function corresponding closely to our function $\phi$ used in
the proof of Theorem~\ref{main}. However, Yoon's approach is simply
to smooth at this stage, obtaining the equivalent of our function
$F$ in the proof of Theorem~\ref{main}. Because there is no
preprocessing of $f$ to $H$, Yoon's function $F$ does not enjoy the
nice property $F(a) = f(a)$ for all $a\in \nodes$. It is this
property which makes the following step, where we treat $G=f-F$, a
fairly simple process. Correspondingly, Yoon has considerably more
difficulty treating his function $G$. Our method also yields the
same bound as that in Yoon, but for a wider class of functions.
Indeed we would suggest that $BL^m(\Omega)$ is the natural class of
functions for which one would wish an error estimate of the type
given in Theorem~\ref{main}.

\end{document}